\documentclass[11pt,psamsfonts,twoside]{article}
\usepackage{amsmath,amssymb,latexsym}
\usepackage{diagrams}
\usepackage{rotating}
\usepackage[mathscr]{eucal}

\newcommand\qed{{\hspace*{\fill}Q.E.D.\vskip12pt plus 1pt}}

\newcommand\Pic[1]{\hbox{\rm Pic(}#1\hbox{\rm )}}

\newcommand\sO{{\mathscr O}}

\def\Coker{\operatorname{Coker}}

\def\rank{\operatorname{rank}}

\def\alg{\operatorname{alg}}

\def\lim{\operatorname{lim}}
\def\Lef{\operatorname{Lef}}
\def\Leff{\operatorname{Leff}}

\def\NS{\operatorname{NS}}
\def\Num{\operatorname{Num}}
\def\Proj{\operatorname{Proj}}

\def\codim{\operatorname{codim}}

\def\Pic{\operatorname{Pic}}

\def\id{\operatorname{id}}
\def\Ker{\operatorname{Ker}}

\def\Coker{\operatorname{Coker}}

\newcommand\pn[1]{{\mathbb P}^{#1}}

\newcommand\proof{\noindent{\em Proof.}\ \ }

\newtheorem{theorem}{Theorem}[section]
\newtheorem{lemma}[theorem]{Lemma}
\newtheorem{corollary}[theorem]{Corollary}

\newtheorem{question}[theorem]{Question}

\newtheorem{definition}[theorem]{Definition}
\newtheorem{rem}[theorem]{Remark}

\newtheorem{pargrph}[theorem]{}
\newtheorem{examp}[theorem]{Example}

\newtheorem{MM}[theorem]{ }
\newtheorem{res}[theorem]{Remarks}

\makeatletter
\ifnum\@ptsize=0  \fi
\ifnum\@ptsize=2  \fi
\renewcommand{\qed}{\hfill $\square$}

\newenvironment{rem*}{\begin{rem}\em}{\end{rem}}
\newenvironment{rems*}{\begin{res}\em}{\end{res}}
\newenvironment{example*}{\begin{examp}\em}{\end{examp}}
\newenvironment{definition*}{\begin{definition}\em}{\end{definition}}
\newenvironment{question*}{\begin{question}\em}{\end{question}}
\newenvironment{MM*}{\begin{MM}\em}{\end{MM}}

\newenvironment{prgrph*}[1]{\indent\begin{pargrph}{\bf #1.}\em\
}{\end{pargrph}}

\begin{document}

\title{A Barth-Lefschetz theorem for submanifolds \\of a product of projective spaces}
\author{Lucian B\u adescu and Flavia Repetto}
\date{}

\maketitle

\begin{abstract} Let $X$ be a complex submanifold of dimension $d$ of $\mathbb P^m\times\mathbb P^n$ ($m\geq n\geq 2$) and 
denote by $\alpha\colon\Pic(\mathbb P^m\times\mathbb P^n)\to \Pic(X)$ the restriction map of Picard groups, by 
$N_{X|\mathbb P^m\times\mathbb P^n}$ the normal bundle of $X$ in $\mathbb P^m\times\mathbb P^n$. Set  
$t:=\max\{\dim\pi_1(X),\dim\pi_2(X)\}$, where $\pi_1$ and $\pi_2$ are the two projections of $\mathbb P^m\times\mathbb P^n$. We prove a Barth-Lefschetz type result as follows:
{\em Theorem.} {\it If $d\geq \frac{m+n+t+1}{2}$ then $X$ is algebraically simply connected, the map $\alpha$ is injective and $\Coker(\alpha)$ is torsion-free. 
Moreover $\alpha$ is an isomorphism if $d\geq\frac{m+n+t+2}{2}$, or if $d=\frac{m+n+t+1}{2}$ and $N_{X|\mathbb P^m\times\mathbb P^n}$ is decomposable.} These bounds are optimal. The main technical ingredients in the proof are: the Kodaira-Le Potier vanishing theorem in the generalized form of Sommese (\cite{LP}, \cite{ShS}), the join construction and an algebraisation result of Faltings concerning small codimensional subvarieties in 
$\mathbb P^N$ (see \cite{Fa}).\end{abstract}

\section*{Introduction} 

It is well known that if $X$ is a submanifold of the complex projective space $\mathbb P^{n}$ 
($n\geq 3$) of dimension $d>\frac{n}{2}$ then a topological result of Lefschetz type, due to Barth and Larsen (see  
\cite{L}, \cite{Ba}), asserts that the canonical restriction maps
$H^i(\mathbb P^{n},\mathbb Z)\to H^i(X,\mathbb Z)$ 
are  isomorphisms for $i\leq 2d-n$, and injective with torsion-free cokernel, for $i=2d-n+1$. As a consequence, the restriction map $\Pic(\mathbb P^{n})\to\Pic(X)$ is an isomorphism if $d\geq\frac{n+2}{2}$, and injective with torsion-free cokernel if $n=2d-1$.

This topological result has been generalized by Sommese to the case when the ambient space $\mathbb P^n$ is replaced by any projective rational homogeneous space $M$ (see \cite{S1}). For example, if $M=\mathbb P^m\times\mathbb P^n$ (with $m\geq n\geq 2$) then Sommese's topological result implies that the canonical restriction map 
$$\alpha\colon\Pic(\mathbb P^m\times\mathbb P^n)\to\Pic(X)$$
is injective with torsion-free cokernel for every submanifold $X$ of $\mathbb P^m\times\mathbb P^n$ of dimension 
$d\geq\frac{2m+n+1}{2}$, and an isomorphism if $d\geq\frac{2m+n+2}{2}$.

The aim of this paper is to prove (in a geometric way) an improved version of Sommese's result concerning the Picard group of the small-codimensional submanifolds 
$X$ of $\mathbb P^m\times\mathbb P^n$ of dimension $d$. To state the main result, set
\begin{equation}\label{t}
t:=\max\{\dim\pi_1(X),\dim\pi_2(X)\},
\end{equation}
where $\pi_1$ and $\pi_2$ are the two canonical projections of $\mathbb P^m\times\mathbb P^n$.

\medskip

\noindent{\bf Main Theorem.} {\em Let $X$ be a complex submanifold of dimension $d$ of $\mathbb P^m\times\mathbb P^n$ 
$($with $m\geq n\geq 2)$, and denote by $\alpha\colon\Pic(\mathbb P^m\times\mathbb P^n)\to \Pic(X)$ the restriction map of Picard groups, by $N_{X|\mathbb P^m\times\mathbb P^n}$ the normal bundle of $X$ in $\mathbb P^m\times\mathbb P^n$, and by $t$ the integer defined by \eqref{t}. Then the following statements hold true:}
\begin{enumerate}
\item[i)] {\it If $d\geq \frac{m+n+t+1}{2}$ then $X$ is algebraically simply connected $($and in particular, 
$H^1(\mathscr O_X)=0)$, the map $\alpha$ is injective and $\Coker(\alpha)$ is torsion-free. }
\item[ii)] {\it  If $d\geq\frac{m+n+t+2}{2}$ then $\alpha$ is an isomorphism.}
\item[iii)] {\it If $d=\frac{m+n+t+1}{2}$ and $N_{X|\mathbb P^m\times\mathbb P^n}\cong E_1\oplus E_2$, with $E_1$ and 
$E_2$ vector bundles of rank $\geq 1$, then $\alpha$ is an isomorphism.}
\item[iv)] {\it If $d=\frac{m+n+t}{2}$ and $N_{X|\mathbb P^m\times\mathbb P^n}\cong E_1\oplus E_2$, with $E_1$ and $E_2$ vector bundles of rank $\geq 1$, then $H^1(\mathscr O_X)=0$. If instead $N_{X|\mathbb P^m\times\mathbb P^n}\cong E_1\oplus E_2$, with $E_1$ and $E_2$ vector bundles of rank $\geq 2$, then $\alpha$ is injective and $\rank\NS(X)=2$, where $\NS(X)$ denotes the N\'eron-Severi group of $X$.}\end{enumerate}

\medskip

We show by examples that the bounds given in Main Theorem, i) and ii) are optimal. As far as i) is concerned, take for instance $m=2s$ even and $s\leq n\leq m=2s$. For every $s\geq 2$ there exist elliptic scrolls $Y$ of dimension $s$ in $\pn m=\pn{2s}$ (well known if $s=2$ and \cite{I} for every $s\geq 3$). Set $X=Y\times\pn n$. Then 
$d=s+n$, $t=n$ and $d=\frac{m+n+t}{2}$. Since $H^1(\mathscr O_X)=H^1(\mathscr O_Y)\neq 0$, $X$ is not algebraically simply connected, and $\Coker(\alpha)$ has a lot of torsion (because 
$\Pic^0(X)\neq 0$). Since $H^1(\mathscr O_X)\neq 0$ the normal bundle $N_{X|\mathbb P^m\times\mathbb P^n}=
q^*(N_{Y|\mathbb P^m})$ (where $q\colon X=Y\times\pn n\to Y$ is the canonical projection) is indecomposable by part iv) of the Main Theorem. However, the fact that $N_{Y|\mathbb P^m}$ is indecomposable was previously proved in \cite{B}. Note also that in this case $\rank\NS(X)=3$. 

To produce an example showing that the bound in ii) is also optimal, take $m=2s+1$ with 
$s\geq 2$ and $n$ such that $s+1\leq n\leq m=2s+1$. Let $Y$ be the image of the Segre embedding $\pn s\times\pn 1\hookrightarrow\pn m=\pn {2s+1}$, and set $X=Y\times\pn n$. Then $d=s+n+1$ and $t=n$ are such that $d=\frac{m+n+t+1}{2}$.
However, the map $\alpha$ cannot be an isomorphism because $\rank\Pic(X)=3$. Notice also that by iii) the normal bundle 
$N_{X|\mathbb P^m\times\mathbb P^n}=q^*(N_{Y|\mathbb P^m})$ is indecomposable (with $q\colon X=Y\times\pn n\to Y$ the canonical projection). The indecomposability of $N_{Y|\mathbb P^m}$ was also previously proved in \cite{B}.
 
The proof of part i) makes use of the join construction to reduce the problem to an open subset of a small-codimensional subvariety of $\mathbb P^{m+n+1}$ and then to apply a result of Faltings (see \cite{Fa}). This is done in section 2. The proof of parts ii)--iv) (which is inspired from \cite{B}) makes systematic use of Kodaira-Le Potier vanishing theorem in the generalized form given by Sommese (see \cite{ShS}) and is contained in section 1.

We want to mention the following interesting recent result of Arrondo and Caravantes \cite{AC} which is related to our Main Theorem (although our approach is completely different from theirs): 

\medskip

\noindent{\bf Theorem} {\bf(Arrondo-Caravantes)} { \em Let $X$ be a complex submanifold of $\pn m\times\pn m$ of dimension 
$d\geq m+1$ such that $\pi_i(X)=\pn m$ for $i=1,2$. Then the restriction map 
$\alpha\colon\Pic(\pn m\times\pn m)\to\Pic(X)$ is injective and $\Pic(X)$ is a free abelian group of rank two.}

\medskip

In the result of Arrondo and Caravantes the codimension of $X$ is relatively larger than in our Main Theorem, but it does not give any information on the torsion of $\Coker(\alpha)$. 

\medskip

All varieties considered throughout are defined over the field $\mathbb C$ of complex numbers. By a manifold we mean a nonsingular irreducible complex algebraic variety. The terminology and the notation used are standard, unless otherwise specified.

\section{A general result on submanifolds in $\mathbb P^m\times\mathbb P^n$}\label{first}\addtocounter{subsection}{1}
\setcounter{theorem}{0}

Let $X$ be a closed irreducible subvariety of dimension $d\geq 1$ of $\mathbb P^m\times\mathbb P^n$, with $m\geq n\geq 
2$. Denote by $\pi_1\colon\mathbb P^m\times\mathbb P^n\to\mathbb P^m$ and $\pi_2\colon\mathbb P^m\times\mathbb P^n\to\mathbb P^n$ the canonical projections and by $p_1\colon X\to p_1(X)$ and $p_2\colon X\to p_2(X)$ the restrictions 
$\pi_1|X$ and $\pi_2|X$. 
Throughout this paper we shall assume that $p_1(X)$ and $p_2(X)$ are both positive dimensional. A closed irreducible subvariety $X$ of $\mathbb P^m\times\mathbb P^n$ satisfying this property will be called {\em positive}. It is well known (and easy to see) that $X$ is positive if and only if $X$ intersects every hypersurface of 
$\mathbb P^m\times\mathbb P^n$. Set
\begin{equation}\label{a,b}
\mathscr O_X(a,b):=\mathscr O_{\mathbb P^m\times\mathbb P^n}(a,b)|X=p_1^*(\mathscr O_{p_1(X)}(a))\otimes p_2^*(\mathscr 
O_{p_2(X)}(b)),\;\forall a,b\in\mathbb Z. 
\end{equation}

\medskip

Since the cotangent bundle of $\mathbb P^m\times\mathbb P^n$ is given by 
$\Omega^1_{\mathbb P^m\times\mathbb P^n}=\pi_1^*(\Omega^1_{\mathbb P^m})\oplus \pi_2^*(\Omega^1_{\mathbb P^n})$,
we get 
\begin{equation}\label{tan}
\Omega^1_{\mathbb P^m\times\mathbb P^n}|X=p_1^*(\Omega^1_{\mathbb P^m}|p_1(X))\oplus p_2^*(\Omega^1_{\mathbb P^n}|p_2(X)). 
\end{equation}

Then a lot of information about the embedding $X\subseteq\mathbb P^m\times\mathbb P^n$ is contained in the following commutative diagram with exact rows and columns:
{\small
\begin{diagram}
&      &&&0&             &0\\
&     &&&\dTo&            &\dTo\\
0&\rTo &N_{X|\mathbb P^{m}\times\mathbb P^n}^{\vee}&\rTo &p_1^*(\Omega^1_{\mathbb P^{m}}|p_1(X))\oplus p_2^*(\Omega^1
_{\mathbb P^{n}}|p_2(X))&\rTo & \Omega^1_X&\rTo &0\\
&      &\dTo^{\id}&              &\dTo&                  &\dTo\\
0&\rTo &N_{X|\mathbb P^{m}\times\mathbb P^n}^{\vee}&\rTo^{\beta'} &p_1^*(\mathscr O_{p_1(X)}(-1)^{\oplus m+1})\oplus p_2^*(\mathscr O_{p_2(X)}(-1)^{\oplus n+1})&\rTo^{\beta}  &F&\rTo &0\\
&      &&&\dTo^{\gamma}&             &\dTo_{\varepsilon}\\
&     && &p_1^*(\mathscr O_{p_1(X)})\oplus p_2^*(\mathscr O_{p_2(X)})=\mathscr O_X^{\oplus 2}&\rTo{\id} &\mathscr O_X^{\oplus 2}\\
&      &&&\dTo&             &\dTo\\
&     && &0&              &0\\
\end{diagram}}
in which $N_{X|\mathbb P^{m}\times\mathbb P^n}^{\vee}$ is the conormal bundle of $X$ in $\mathbb P^{m}\times
\mathbb P^n$, the first row is the canonical exact sequence of cotangent bundles of $X$ in $\mathbb P^m\times\mathbb P^n$ (taking into account of \eqref{tan}), the middle column is the direct sum of the restricted Euler sequences of $\mathbb P^m$ and of $\mathbb P^n$, and $F:=\Coker(\beta')$.

We first prove the following general result:

\begin{lemma}\label{gen1.5} Let $f\colon X\to Y$ be a proper surjective morphism from a projective manifold $X$ of dimension $d\geq 2$ onto a projective variety $Y$ of dimension $e$ with $d>e\geq 2$. Then $f$ has no fiber isomorphic to 
$\mathbb P^{d-1}$.\end{lemma}

\proof Assume that there exists $y\in Y$ such that the fiber $F=f^{-1}(y)$ is isomorphic to $\mathbb P^{d-1}$. Then $F$ is an effective divisor on $X$. Since $\Pic(\mathbb P^{d-1})=\mathbb Z[\mathscr O_{\mathbb P^{d-1}}(1)]$, it follows that the conormal line bundle $N_{F|X}^{\vee}$ is isomorphic to $\mathscr O_{\mathbb P^{d-1}}(s)$ for some $s\in\mathbb Z$.

We claim that $s>0$. Indeed, since $\Pic(\mathbb P^{d-1})\cong\mathbb Z$, it is sufficient to show that we can find an irreducible curve $C\subseteq F$ such that $\deg_C(N_{F|X}^{\vee}|C)=(N_{F|X}^{\vee}\cdot C)>0$. 
To produce such a curve (following an idea of P. Ionescu) we fix a projective embedding $X\hookrightarrow\mathbb P^N$, and let $H_1,\ldots,H_{d-2}$ be $d-2$ general hyperplanes of $\mathbb P^N$, and set $X':=X\cap H_1\cap\cdots\cap H_{d-2}$ and $C:=F\cap H_1\cap\cdots\cap H_{d-2}$. By Bertini, $X'$ is a smooth projective surface and $C$ is a smooth irreducible curve on $X'$. By construction, the morphism $f':=f|X'\colon X'\to f(X')$ is generically finite and $f'(C)=y$; then by a well known elementary fact in the theory of surfaces, $(C^2)_{X'}<0$ (for instance this fact is an easy consequence of Hodge index theorem). On the other hand, since $C$ is the proper intersection of $F$ with $H_1\cap\cdots\cap H_{d-2}$ we infer that $N_{F|X}|X'\cong N_{C|X'}$. But since $(C^2)_{X'}=\deg_C(N_{C|X'})=
\deg_C(N_{F|X})$, we get $\deg_C (N_{F|X}^{\vee})|C>0$, as claimed.

Now, by a generalization of a contractibility result of Castelnuovo-Kodaira (see \cite{B3}) the divisor $F=\mathbb P^{d-1}$ of $X$ of conormal bundle $\mathscr O_{\mathbb P^{d-1}}(s)$ with $s>0$ can be blown down to a normal point, i.e. there exists a birational morphism $\varphi\colon X\to V$, with $V$ a normal projective variety such that $\varphi(F)=v$ is a point and $\varphi|X\setminus F$ defines a biregular isomorphism $X\setminus F\cong V\setminus
\{v\}$. Then by a well known elementary fact, there is a unique morphism $g\colon V\to Y$ such that $g\circ\varphi=f$. In particular, the fiber $g^{-1}(g(v))$ is reduced to the point $v$, which contradicts the theorem on the dimension of fibers
because by hypothesis $\dim(V)=d>e=\dim(Y)$.\qed 

\begin{corollary}\label{gen1.6} Let $X$ be a submanifold of dimension $d\geq m+1$ of $\mathbb P^{m}\times\mathbb 
P^n$, with $m\geq n\geq 2$. If $d=m+1$ assume moreover that $\dim p_i(X)\geq 2$ for $i=1,2$. Then all the fibers of 
$p_i\colon X\to p_i(X)$, $i=1,2$, have dimension $\leq d-2$.\end{corollary}

\proof The assertion is trivial if $d\geq m+2$ because all the fibers of $p_i$ ($i=1,2$) are of dimension $\leq m$. Assume therefore $d=m+1$; if there exists a fiber $F$ of $p_i$ of dimension $m$ (with $i=1$ or $i=2$), then necessarily $F\cong\mathbb P^m$. In this case the corollary follows from Lemma \ref{gen1.5}.\qed

\medskip

Now we turn to our general situation (under the hypotheses from the beginning). By \eqref{a,b} we have 
$$p_1^*(\mathscr O_{p_1(X)}(-1)^{\oplus m+1})\cong\mathscr O_X(-1,0)^{\oplus m+1}\;\;\text{and}\;\;
p_2^*(\mathscr O_{p_2(X)}(-1)^{\oplus n+1})\cong\mathscr O_X(0,-1)^{\oplus n+1}.$$ 
Moreover,
$$H^0(p_1^*(\mathscr O_{p_1(X)}(-1)^{\oplus m+1})\oplus p_2^*(\mathscr O_{p_2(X)}(-1)^{\oplus n+1}))=0.$$ 
This follows because $p_1(X)$ and $p_2(X)$ are positive dimensional, whence by the above isomorphisms 
$\mathscr O_X(1,0)^{\oplus m+1}$ and $\mathscr O_X(0,1)^{\oplus n+1}$ are direct sums of $(d-1)$-ample line 
bundles (in the sense of Sommese \cite{S}). 
Thus the cohomology of the second row of the above diagram yields the exact sequence
\begin{equation}\label{e1}\begin{split}
&0\to H^0(F)\to H^1(N_{X|\mathbb P^{m}\times\mathbb P^n}^{\vee})\to 
H^1(\mathscr O_X(-1,0)^{\oplus m+1}\oplus\mathscr O_X(0,-1)^{\oplus n+1})\to \\
& \to H^1(F)\to H^2(N_{X|\mathbb P^{m}\times\mathbb P^n}^{\vee})\to H^2(\mathscr O_X(-1,0)^{\oplus m+1}\oplus\mathscr 
O_X(0,-1)^{\oplus n+1}).\end{split}\end{equation}

On the other hand, the last column yields the cohomology exact sequence
\begin{equation}\label{e2}
0\to H^0(\Omega^1_X)\to H^0(F)\to H^0(\mathscr O_X^{\oplus 2})\to H^1(\Omega^1_X)\to H^1(F)\to 
H^1(\mathscr O_X^{\oplus 2}).\end{equation}

\begin{lemma}\label{van1} Under the above hypotheses, assume moreover that the projections $p_i\colon X\to p_i(X)$ have all fibers of dimension $\leq d-2$ for $i=1,2$, e.g. if $d\geq m+2$, or if $d=m+1$ and  $\dim p_i(X)\geq 2$ for $i=1,2$, 
$($by Corollary $\ref{gen1.6}$ above$)$. Then 
$H^1(\mathscr O_X(-1,0)^{\oplus m+1}\oplus\mathscr O_X(0,-1)^{\oplus n+1})=0$.\end{lemma}

\proof The hypothesis implies that $\mathscr O_X(1,0)^{\oplus m+1}\cong p_1^*(\mathscr O_{p_1(X)}(1)))^{\oplus m+1}$ and 
$\mathscr O_X(0,1)^{\oplus n+1}\cong p_2^*(\mathscr O_{p_2(X)}(1))^{\oplus n+1}$ are both $(d-2)$-ample vector bundles which are direct sums of line bundles. Then the conclusion follows from Kodaira vanishing theorem in the generalized form of Sommese, see \cite{ShS}, page 96, Corollary (5.20).\qed

\begin{corollary}\label{van2} Under the hypotheses of Lemma $\ref{van1}$ one has $h^1(N_{X|\mathbb P^{m}\times
\mathbb P^n}^{\vee})= h^0(F)$ and $h^1(F)\leq h^2(N_{X|\mathbb P^{m}\times\mathbb P^n}^{\vee})$.
\end{corollary}

\proof The corollary follows from the exact sequence \eqref{e1} and from Lemma \ref{van1}.\qed

\begin{corollary}\label{van3} Under the hypotheses of Lemma $\ref{van1}$ one has $h^1(\Omega^1_{\mathbb P^m\times
\mathbb P^n}|X)=2$.\end{corollary}

\proof The corollary follows from the cohomology sequence of the first column of the above diagram  and from
Lemma \ref{van1}, taking into account of the isomorphism \eqref{tan}.\qed

\medskip

Now we analyze the exact sequence \eqref{e2}. The pull backs of the Euler exact sequences
{\small
\begin{diagram}
0&\rTo &p_1^*(\Omega^1_{\mathbb P^m}|p_1(X))&\rTo &p_1^*(\mathscr O_{p_1(X)}(-1)^{\oplus m+1})=\mathscr O_X(-1,0)^{\oplus m+1}&\rTo^{\gamma_1}&\mathscr O_X&\rTo& 0,\end{diagram}}
{\small
\begin{diagram}
0&\rTo &p_2^*(\Omega^1_{\mathbb P^n}|p_2(X))&\rTo &p_2^*(\mathscr O_{p_2(X)}(-1)^{\oplus n+1})=\mathscr O_X(0,-1)^{\oplus n+1}&\rTo^{\gamma_2}&\mathscr O_X&\rTo& 0\end{diagram}}
\medskip
do not split, because $p_1(X)$ and $p_2(X)$ are positive dimensional. This means that $H^0(\gamma_i)=0$ for $i=1,2$. Since the first vertical column of the above diagram is the direct sum of these exact sequences, it follows that the map
$$H^0(\gamma)=H^0(\gamma_1)\oplus H^0(\gamma_2)$$
is also zero. Thus from the cohomology sequence of the first column we infer that the map 
$\delta_1\colon H^0(\mathscr O_X^{\oplus 2})\to H^1(p_1^*(\Omega^1_{\mathbb P^m}|p_1(X))\oplus 
p_2^*(\Omega^1_{\mathbb P^n}|p_2(X)))$ is injective.

On the other hand, in the commutative square
{\small
\begin{diagram}
H^0(\mathscr O_X^{\oplus 2})&\rTo^{\id}&H^0(\mathscr O_X^{\oplus 2})\\
\dTo^{\delta_1}&&\dTo_{\delta_2}\\
H^1(p_1^*(\Omega^1_{\mathbb P^m}|p_1(X))\oplus p_2^*(\Omega^1_{\mathbb P^n}|p_2(X)))&\rTo&H^1(\Omega^1_X)\\
\end{diagram}}
the bottom horizontal map is not zero. Indeed, by hypothesis $p_1(X)$ and $p_2(X)$ are both positive dimensional. Since 
$\dim p_1(X)>0$ the map $H^1(\Omega^1_{\mathbb P^m})\to H^1(\Omega^1_X)$ is non-zero (the image of the class of 
$\mathscr O_{\mathbb P^m}(1)$ is not zero in $H^1(\Omega^1_X)$). Since this map is the composition
$$H^1(\Omega^1_{\mathbb P^m})\to H^1(p_1^*(\Omega^1_{\mathbb P^m}|p_1(X)))\to H^1(\Omega^1_X),$$ 
it follows that the second map cannot be zero.

Now, since $\delta_1$ is injective we infer that $\delta_2\neq 0$. Thus \eqref{e2} yields the exact sequences
\begin{equation}\label{irr}0\to H^0(\Omega^1_X)\to H^0(F) \;\;\text{and}\;\;
0\to V\to H^1(\Omega^1_X)\to H^1(F)\to H^1(\mathscr O_X^{\oplus 2}),\end{equation}
in which $V$ is a $\mathbb C$-vector space of dimension $2$ if $\delta_2$ is injective (i.e. if $H^0(\varepsilon)=0$), and $1$ otherwise. In particular,
\begin{equation}\label{e3}
h^1(\Omega^1_X)\leq\begin{cases}2+h^1(F),\;\;\text{if}\;\; H^0(\varepsilon)=0\\
1+h^1(F),\;\;\text{if}\;\; H^0(\varepsilon)\neq 0.\end{cases}
\end{equation}
Moreover, in both cases we have equality if $H^1(\mathscr O_X)=0$.

Putting everything together and using Corollary \ref{van2} we get:

\begin{theorem}\label{gen} Under the hypotheses of the beginning assume moreover that both projections 
$p_i\colon X\to p_i(X)$, $i=1,2$, have fibers all of dimension $\leq d-2$. $($This is always the case if $d\geq m+2$, or by Corollary $\ref{gen1.6}$ above, if $d=m+1$ and $\dim p_i(X)\geq 2$ for $i=1,2.)$ Then the following statements hold true:
\begin{enumerate}
\item[\em i)] $h^{0}(\Omega^1_X)\leq h^1(N_{X|\mathbb P^{m}\times\mathbb P^n}^{\vee})$. 
\item[\em ii)]$\rank\NS(X)\leq\begin{cases}2+h^2(N_{X|\mathbb P^{m}\times\mathbb P^n}^{\vee}),\;\;\text{if}\;\; H^0(\varepsilon)=0,\\
1+h^2(N_{X|\mathbb P^{m}\times\mathbb P^n}^{\vee}),\;\;\text{if}\;\; H^0(\varepsilon)\neq 0,\end{cases}$

where $\NS(X)$ is the N\'eron-Severi group of $X$. 
\item[\em iii)] If $d\geq m+1$ and $H^1(N_{X|\mathbb P^{m}\times\mathbb P^n}^{\vee})=0$ then $H^1(\mathscr O_X)=0$ and 
$2\leq\rank\NS(X)\leq 2+h^2(N_{X|\mathbb P^{m}\times\mathbb P^n}^{\vee})$.\end{enumerate}\end{theorem}

\proof i) follows from Corollary \ref{van2} and from \eqref{irr}.
ii) follows from the well known inequality $\rank\NS(X)\leq h^1(\Omega^1_X)$ (valid in characteristic zero, see 
\cite{H}, Exercise 1.8, page 367, if $X$ is a surface, and \cite{B}, the claim in the proof of Theorem 2.1, in general), and from \eqref{e3}.  The assertion about $H^1(\mathscr O_X)$ in iii) follows from i), using the Hodge symmetry 
$H^0(\Omega^1_X)\cong H^1(\mathscr O_X)$ (via Serre's GAGA). The last part of iii) follows from Lemma \ref{inj}, and from ii) because $H^0(\varepsilon)=0$ if $H^1(N_{X|\mathbb P^{m}\times\mathbb P^n}^{\vee})=0$, by Corollary \ref{van2}.\qed

\medskip

\begin{rem*}\label{gen1} Take $m=n\geq 2$ and $X=\Delta\cong\mathbb P^n$ the diagonal of $\mathbb P^{n}\times
\mathbb P^n$.  Then $N_{X|\mathbb P^{n}\times\mathbb P^n}^{\vee}=\Omega^1_{\mathbb P^n}$, whence
$H^1(N_{X|\mathbb P^{n}\times\mathbb P^n}^{\vee})=H^1(\Omega^1_{\mathbb P^n})\cong\mathbb C$.
So by Theorem \ref{gen}, i), $H^0(F)\cong\mathbb C$. In this case,  $H^0(\varepsilon)\neq 0$, whence
$1=\rank\Num(X)\leq 1+h^1(F)$.
Moreover, $H^2(N_{X|\mathbb P^{n}\times\mathbb P^n}^{\vee})=H^2(\Omega^1_{\mathbb P^n})=0$, whence $h^1(F)=0$ by Corollary
\ref{van2}. Thus the above inequality becomes equality. In this case we also have $H^0(\Omega^1_X)=0$ and 
$H^1(N_{X|\mathbb P^{n}\times\mathbb P^n}^{\vee})\neq 0$. In particular, the second possibility in \eqref{e3} really  occurs.\end{rem*}

\begin{corollary}\label{gen3} Assume that $N_{X|\mathbb P^{m}\times\mathbb P^n}$ is ample and $d\geq m+1$. Then the irregularity of $X$ is zero, and  $\rank\Pic(X)=2$.\end{corollary}

\proof Since $N_{X|\mathbb P^{m}\times\mathbb P^n}$ is ample from Le Potier vanishing theorem and $d\geq m+1$ it follows that $H^i(N_{X|\mathbb P^{m}\times\mathbb P^n}^{\vee})=0$ for $i\leq 2$.
Then the conclusion follows from Theorem \ref{gen}.\qed

\medskip

In the sequel we shall be interested in the submanifolds of $\mathbb P^{m}\times\mathbb P^n$ of dimension $d$ with
$d\geq\frac{m+n+t+1}{2}$. We shall need the following lemma:

\begin{lemma}\label{gen5} Let $X$ be a submanifold of $\mathbb P^{m}\times\mathbb P^n$ of dimension $d$ with 
$m\geq n\geq 2$. 
\begin{enumerate}
\item[\em i)] Assume that $d\geq\frac{m+n+t+1}{2}$. Then $d\geq m+1$ and if $d=m+1$ then $\dim p_i(X)\geq 2$, $i=1,2$.
\item[\em ii)] Assume that $d=\frac{m+n+t}{2}$. Then $d\geq m$. If $d=m$ then $m\geq 2n$ and $X=X_1\times\mathbb P^n$,
with $X_1$ a submanifold of $\mathbb P^m$ of dimension $m-n$, with $m-n\geq n\geq 2$. If $d\geq m+1$ and $\dim p_i(X)=1$ for some $i\in\{1,2\}$ then $X$ is isomorphic to $\mathbb P^m\times X_2$, with $X_2$ a smooth curve in $\mathbb P^n$, with $n=2$.\end{enumerate}\end{lemma}

\proof i) Assume $d\leq m$, i.e. $d=m-r$, with $r\geq 0$. The hypothesis implies $m+n+t+1\leq 2d=2m-2r$, whence 
$d-t=m-r-t\geq n+r+1$. Let $F_1$ be a general fiber of $p_1\colon X\to p_1(X)$. Since $t\geq\dim p_1(X)$, by the theorem on dimension of fibers we get
$$\dim(F_1)=d-\dim p_1(X)\geq d-t\geq n+r+1\geq n+1.$$
However this is impossible because $F_1\cong p_2(F_1)\subseteq\mathbb P^n$. This proves that $d\geq m+1$.

Now we prove the last part of i). The only case in which we can have $\dim p_i(X)=1$ for some $i\in\{1,2\}$ is when 
$d=m+1$. Thus 
\begin{equation}\label{w} d=m+1\geq \frac{m+n+t+1}{2}.\end{equation}
If $m=n$ the inequality \eqref{w} implies $t=1$. On the other hand, since $X\subseteq p_1(X)\times p_2(X)$, it follows that $d\leq 2$,
and in particular, by the above inequality we get $m=n=1$, which contradicts the hypothesis that $m\geq n\geq 2$. If instead $m>n$ then $\dim(F_1)\leq n$, whence
$\dim p_1(X)=m+1-\dim(F_1)\geq m+1-n\geq 2$. Assume now that $\dim p_2(X)=1$; then all fibers of $p_2$ are 
$m$-dimensional, whence all of them are isomorphic to $\mathbb P^m$ (since they are contained in $\mathbb P^m\times p$, with $p\in\mathbb P^n$). It follows that $X=p_2^{-1}(p_2(X))=\mathbb P^m\times
p_2(X)$ and in particular, $t=m$. Finally, since $m+1\geq \frac{m+n+t+1}{2}=m+\frac{n+1}{2}$ we get $n=1$, which again contradicts the hypotheses. This proves i).

\medskip

ii) Assume first that $d\leq m$, i.e. $d=m-r$, with $r\geq 0$. Since $t\geq\dim p_1(X)$ and $F_1\cong p_2(F_1)\subseteq\mathbb P^n$, the equality $d=\frac{m+n+t}{2}$ and the theorem on dimension of fibers yield 
 $$n+r=d-t\leq d-\dim p_1(X)=\dim(F_1)\leq n,$$
where (as above) $F_1$ is a general fiber of $p_1\colon X\to p_1(X)$. It follows that $r=0$, i.e. 
$d=m$, $t=\dim p_1(X)$ and $F_1=\mathbb P^n$. Hence all fibers of $p_1$ are isomorphic to $\mathbb P^n$, i.e. $X=X_1\times\mathbb P^n$, with $X_1=p_1(X)$ a submanifold of $\mathbb P^m$ of dimension $m-n$. Moreover, since 
$t=\dim p_1(X)=d-n=m-n$ and $\dim p_2(X)=n$ it follows that $m-n\geq n$, i.e. $m\geq 2n$. In particular, every fiber of the projections $p_1$ and $p_2$ is of dimension $\leq d-2$.

Assume now $d\geq m+1$. Then as above we can have $\dim p_i(X)=1$ for some $i\in\{1,2\}$ only if $d=m+1$. Thus 
$m+1=\frac{m+n+t}{2}$ yields $m=n+t-2$. If $m=n$ then $t=2$, and therefore (say) $\dim p_1(X)=2$ and $\dim p_2(X)=1$.
Using $X\subseteq p_1(X)\times p_2(X)$ this immediately yields $d\leq 3$ and $X$ is the hypersurface $\mathbb P^2\times p_2(X)\subset\mathbb P^2\times\mathbb P^2$. If $m>n$ then only $X_2:=p_2(X)$ can be a curve, and in this case $X=\mathbb P^m\times X_2$. Since in this case $t=m$ it follows $m=n+m-2$, i.e. $n=2$ and $m\geq 3$. \qed

\medskip

\medskip

Now we come back to the above commutative diagram with exact rows and columns. Then from the second row of this diagram it follows that 
$N_{X|\mathbb P^{m}\times\mathbb P^n}$ is a quotient of 
\begin{equation}\label{vb}
p_1^*(\mathscr O_{p_1(X)}(1)^{\oplus m+1})\oplus p_2^*(\mathscr O_{p_2(X)}(1)^{\oplus n+1})=
\mathscr O_X(1,0)^{\oplus m+1}\oplus\mathscr O_X(0,1)^{\oplus n+1}.\end{equation}

Clearly $t\leq m$, where $t$ is defined by formula \eqref{t} of the introduction. Moreover, it is easy to see that the fibers of the morphisms $p_1\colon X\to p_1(X)$ and
$p_2\colon X\to p_2(X)$ are all of dimension $\leq t$. Indeed, if for example $F$ is a fiber of $p_1$ then $F\cong 
p_2(F)$ (via $p_2$), whence $\dim(F)=\dim p_2(F)\leq\dim p_2(X)\leq t$. 

It follows that the vector bundle \eqref{vb} is $t$-ample, whence its quotient $N_{X|\mathbb P^{m}\times\mathbb P^n}$ is a $t$-ample vector bundle of rank $m+n-d$. 
Then using Le Potier-Sommese vanishing theorem (see \cite{ShS}, page 96, Corollary (5.20)) we get:
$$H^i(N_{X|\mathbb P^{m}\times\mathbb P^n}^{\vee})=0,\;\;\text{for}\;\;i\leq d-(m+n-d)-t =2d-m-n-t. $$
In particular, 
{\small
\begin{equation}\label{lps1}\begin{split}
H^1(N_{X|\mathbb P^{m}\times\mathbb P^n}^{\vee})=0\;\;\text{if}\;\;d\geq\frac{m+n+t+1}{2}\;\;\text{and}\\
H^2(N_{X|\mathbb P^{m}\times\mathbb P^n}^{\vee})=0\;\;\text{if}\;\;d\geq\frac{m+n+t+2}{2}.
\end{split}\end{equation}}
Now we are ready to prove the following:

\begin{theorem}\label{gen4}  Let $X$ be a submanifold of $\mathbb P^{m}\times\mathbb P^n$ of dimension $d$ with $m\geq n\geq 2$. Then the following statements hold true:
\begin{enumerate}
\item[\em i)] If $d\geq\frac{m+n+t+1}{2}$ then $H^1(\mathscr O_X)=0$, and if $d\geq\frac{m+n+t+2}{2}$ then 
$\rank\Pic(X)=2$.
\item[\em ii)] If $d=\frac{m+n+t+1}{2}$ and $N_{X|\mathbb P^m\times\mathbb P^n}\cong E_1\oplus E_2$, with $E_1$ and $E_2$ vector bundles of rank $\geq 1$, then $\rank\Pic(X)=2$.
\item[\em iii)] If $d=\frac{m+n+t}{2}$ and  $N_{X|\mathbb P^m\times\mathbb P^n}\cong E_1\oplus E_2$, with $E_1$ and $E_2$ vector bundles of rank $\geq 1$, then $H^1(\mathscr O_X)=0$. If moreover $N_{X|\mathbb P^m\times\mathbb P^n}\cong E_1\oplus E_2$, with $E_1$ and $E_2$ vector bundles of rank $\geq 2$, then $\rank\Pic(X)=2$.\end{enumerate}\end{theorem}

\proof Assume first that $d\geq\frac{m+n+t+1}{2}$. Then by Lemma \ref{gen5}, i) we have $d\geq m+1$ and 
$\dim p_i(X)\geq 2$ for $i=1,2$. Then by Serre's GAGA and Hodge symmetry, $h^0(\Omega^1_X)=h^1(\mathscr O_X)$. Moreover, 
$\NS(X)=\Pic(X)$ if $X$ is a regular variety. Then i) follows from \eqref{lps1} and from Theorem \ref{gen}, i) and iii).

\medskip

ii) By hypothesis $N_{X|\mathbb P^{m}\times\mathbb P^n}=E_1\oplus E_2$, with $E_1$ and $E_2$ vector bundles of rank 
$\leq\rank(N_{X|\mathbb P^{m}\times\mathbb P^n})-1=m+n-d-1$. Since $N_{X|\mathbb P^{m}\times\mathbb P^n}$ is $t$-ample, 
$E_1$ and $E_2$ are also $t$-ample. Thus by Le Potier vanishing theorem in the generalized form given by Sommese (see \cite{ShS}, page 96, Corollary (5.20)), we have
$$H^2(N_{X|\mathbb P^{m}\times\mathbb P^n}^{\vee})\cong H^2(E_1^{\vee})\oplus H^2(E_2^{\vee})=0,$$
because in this case $d-\rank(E_i)-t\geq d-(m+n-d-1)-t=2d-(m+n+t-1)=(m+n+t+1)-(m+n+t-1)=2$, for $i=1,2$. Then by Theorem 
\ref{gen}, iii), $\rank\Pic(X)=2$.

\medskip

iii) By Lemma \ref{gen5}, ii) we may assume that both projections $p_i\colon X\to p_i(X)$, $i=1,2$, have fibers all of dimension $\leq d-2$. Indeed, if $d=m$ by Lemma \ref{gen5}, ii) we have $X=X_1\times\mathbb P^n$ with 
$\dim(X_1)=m-n\geq n\geq 2$. If instead $d\geq m+1$ then we can have $\dim p_i(X)=1$ for some $i\in\{1,2\}$ only if 
$X=\mathbb P^m\times X_2$, with $X_2$ a smooth curve in $\mathbb P^n$, with $n=2$, in which case $X$ is a hypersurface in $\mathbb P^m\times\mathbb P^n$. However this situation is ruled out by the hypotheses which imply 
$\codim_{\mathbb P^m\times\mathbb P^n}(X)\geq 2$.

If $d=\frac{m+n+t}{2}$ and $N_{X|\mathbb P^m\times\mathbb P^n}=E_1\oplus E_2$, with $E_1$ and $E_2$ vector bundles of rank $\leq\rank(N_{X|\mathbb P^{m}\times\mathbb P^n})-1=m+n-d-1$, then as above,
$$H^1(N_{X|\mathbb P^{m}\times\mathbb P^n}^{\vee})\cong H^1(E_1^{\vee})\oplus H^1(E_2^{\vee})=0,$$ 
by Le Potier-Sommese vanishing theorem, because $d=\frac{m+n+t}{2}$ implies $d-\rank(E_i)-t\geq 1$, for $i=1,2$. Then the statement follows from Theorem \ref{gen}, i).

If instead $N_{X|\mathbb P^m\times\mathbb P^n}=E_1\oplus E_2$, with $E_1$ and $E_2$ vector bundles of rank $\geq 2$, i.e. of rank $\leq\rank(N_{X|\mathbb P^{m}\times\mathbb P^n})-2=m+n-d-2$. Then as in the first part of iii), by Le Potier-Sommese vanishing theorem we have 
$$H^2(N_{X|\mathbb P^{m}\times\mathbb P^n}^{\vee})\cong H^2(E_1^{\vee})\oplus H^2(E_2^{\vee})=0,$$ 
because $d=\frac{m+n+t}{2}$ implies $d-\rank(E_i)-t\geq 2$, for $i=1,2$. Then the statement follows from Theorem 
\ref{gen}, ii).\qed

\section{Torsion-freeness of $\Coker(\alpha)$}\label{second}\addtocounter{subsection}{1}
\setcounter{theorem}{0}

In this section we shall prove the following:

\begin{theorem}\label{coker} Let $X$ be a submanifold of dimension $d$ of $\pn m\times\pn n$ $($with $m\geq n\geq 2)$.  If 
$d\geq\frac{m+n+t+1}{2}$ then $X$ is algebraically simply connected, the restriction map $\alpha\colon\Pic(\mathbb P^m\times\mathbb P^n)\to\Pic(X)$ is injective and $\Coker(\alpha)$ is torsion-free, where $t$ is defined by formula \eqref{t} of the introduction. \end{theorem}

We start with the following simple observation:

\begin{lemma}\label{inj} Let $X$ be a positive closed irreducible subvariety of $\mathbb P^m\times\mathbb P^n$ 
$(m\geq n\geq 1)$ such that at least one of the morphisms $p_i\colon X\to p_i(X)$, $i=1,2$ has a positive-dimensional fiber, e.g. if $d>n$. Then the restriction map $\alpha\colon\Pic(\mathbb P^m\times\mathbb P^n)\to\Pic(X)$ is injective.
\end{lemma}

\proof Assume for instance that $p_1$ has a positive-dimensional fiber $F$; by hypotheses we also have $\dim p_1(X)>0$. Let $(a,b)\in\mathbb Z\times\mathbb Z$ such that $\mathscr O_X(a,b)\cong \mathscr O_X$. 
Then this isomorphism implies $\mathscr O_X(a,b)|F\cong\mathscr O_{F}$, and since $\mathscr O_X(a,b)|F \cong \mathscr
O_{F}(a,b)\cong\mathscr O_F(b)$, we get $b=0$ because $\dim(F)>0$ and the restriction map $p_2|F\colon F\to p_2(F)$ is an isomorphism. Thus $\mathscr O_X(a,b)=p_1^*(\mathscr O_{p_1(X)}(a))$. Finally, from $p_1^*(\mathscr O_{p_1(X)}(a))\cong \mathscr O_X$ and $\dim p_1(X)>0$ we get $a=0$ because if $a\neq 0$ one of the line bundles $p_1^*(\mathscr O_{p_1(X)}(a))$ or $p_1^*(\mathscr O_{p_1(X)}(-a))$ is $(d-1)$-ample, while $\mathscr O_X$ is not.\qed

\begin{rem*}\label{inj1} In Lemma \ref{inj} the hypothesis that one of  $p_1$ or $p_2$ has a positive-dimensional fiber is essential; indeed, if we take $m=n$ and $X$ the diagonal of $\mathbb P^m\times\mathbb P^m$ then 
$\Pic(X)\cong\mathbb Z$, while $\Pic(\mathbb P^m\times\mathbb P^m)=\mathbb Z\times\mathbb Z$.\end{rem*}

\begin{corollary}\label{inj'} If $X$ is a closed irreducible subvariety of $\mathbb P^m\times\mathbb P^n$ of dimension
$d\geq\frac{m+n+t}{2}$ $($with $m\geq n\geq 2$ and $t$ given by formula \eqref{t}$)$ then the restriction map $\alpha\colon\Pic(\mathbb P^m\times\mathbb P^n)\to\Pic(X)$ is injective.
\end{corollary}

\proof We shall show that the hypotheses of Lemma \ref{inj} are fulfilled. First we observe that the hypotheses that 
$d\geq\frac{m+n+t}{2}$ and $m\geq n\geq 2$ imply $d>n$. Therefore by the theorem on the dimension of fibers we get 
$\dim p_1(X)>0$ and the morphism $p_2\colon X\to p_2(X)$ has all fibers positive-dimensional. Thus it remains to show that $\dim p_2(X)>0$.  Assuming $\dim (p_2(X))=0$, i.e. $X\subseteq \pn m\times\{p\}\cong\pn m$, with $p\in \pn n$, then $X\cong p_1(X)$, and in particular, $t=d\leq m$. Then the hypothesis that $d\geq\frac{m+n+t}{2}$ yields $d\geq m+n$, a contradiction.\qed

\medskip

Now, to prove the non-trivial part of Theorem \ref{coker} we need some preparation.
Let $X$ be a submanifold of $\pn m\times\pn n$ $($with $m\geq n\geq1)$.
Let us recall the join construction of $X$. Note that this construction has been already used in algebraic geometry in various circumstances, e.g. by Lascu and Scott in \cite{LS} to determine the behaviour of Chern classes undergoing a blowing up, by Deligne in \cite{De} to simplify Fulton-Hansen connectedness theorem \cite{FH}, and by the first named author in \cite{B94} to prove Lefschetz-type results for proper intersections.
In the projective space $\pn {m+n+1}:=\Proj (k[x_0,\dots,x_m,y_0,\dots,y_n])$ consider the disjoint linear subspaces
$$L_1:=\{[x_0,\dots,x_m,y_0,\dots,y_n]\in\pn {m+n+1} |x_0=\dots=x_m=0\},$$
$$L_2:=\{[x_0,\dots,x_m,y_0,\dots,y_n]\in\pn {m+n+1} |y_0=\dots=y_n=0\},$$
and set $U:=\pn {m+n+1}\setminus (L_1\sqcup L_2)$. Consider also the rational map 
$$\pi \colon\pn {m+n+1}\dashrightarrow\pn m\times\pn n$$
defined by
$$\pi([x_0,\dots,x_m,y_0,\dots,y_n]):=([x_0,\dots,x_m],[y_0,\dots,y_n]).$$
Then $\pi$ is defined precisely on $U$ and it is the projection of a locally trivial $\mathbb{G}_m$-bundle (in the Zariski topology), where $\mathbb{G}_m$ is the multiplicative group of $k$. Observe also that the rational map
$\pi_{L_1}:=\pi_1\circ\pi\colon\mathbb P^{m+n+1}\dasharrow\mathbb P^m$ (resp. $\pi_{L_2}:=\pi_2\circ\pi\colon\mathbb P^{m+n+1}\dasharrow\mathbb P^n$) is nothing but the linear projection  of $\mathbb P^{m+n+1}$  of center $L_1$ (resp. the linear projection  of $\mathbb P^{m+n+1}$  of center $L_2$). In particular, $\pi_{L_2}|L_1$ defines an isomorphism $L_1\cong\mathbb P^n$ and $\pi_{L_1}|L_2$ an isomorphism $L_2\cong\mathbb P^m$. Moreover,
$$\pi^*(\pi_1^*(\mathscr O_{\mathbb P^m}(1)))=\mathscr O_U(1)\;\;\text{and}\;\;\pi^*(\pi_2^*(\mathscr O_{\mathbb P^n}(1)))=\mathscr O_U(1).$$

In particular, for every closed irreducible subvariety $X$ of $\mathbb P^m\times\mathbb P^n$ one has
\begin{equation}\label{pr}\pi_X^*(\mathscr O_{X}(1,0))=\mathscr O_{U_X}(1)\;\;\text{and}\;\;
\pi_X^*(\mathscr O_{X}(0,1)))=\mathscr O_{U_X}(1),\end{equation}
where $U_X:={\pi}^{-1}(X)$ and  $\pi_X\colon U_X\to X$ the restriction of $\pi$. Since $\pi\colon U\to\mathbb P^m\times\mathbb P^n$ is a locally trivial $\mathbb G_m$-bundle, so is $\pi_X\colon U_X\to X$. In particular, $U_X$ is irreducible. Denote by $Y:=\overline{U_X}$ the closure of $U_X$  in $\pn {m+n+1}$, then 
$$Y=U_X\sqcup Z,$$
with $Z=(Y\cap L_1)\sqcup(Y\cap L_2)$ (disjoint union).
The following well known fact follows easily:

\begin{lemma}\label{proj} Let $X$ be a closed irreducible subvariety of $\pn m\times\pn n$ $($with $m\geq n\geq 2)$ of dimension $d>0$. In the above notation, one has $Y\cap L_1\cong p_2(X) \ and\ Y\cap L_2\cong p_1(X).$\end{lemma}

We also have the following well known lemma (which is very similar to Lemma 3 in \cite{B94}):

\begin{lemma}\label{L3} Under the above notation, let $P_X:=\mathbb{P}(\sO_X(1,0)\oplus\sO_X(0,1))$ be the projective bundle associated to $\sO_X(1,0)\oplus\sO_X(0,1)$, and denote by $p_X\colon P_X\to X$ the canonical projection of $P_X$. Then the variety $U_X$ can be canonically embedded  in $P_X$ as an open dense subset such that the morphism $\pi_X\colon 
U_X\to X$ extends to $p_X\colon P_X\to X$, and the complement of $U_X$ in $P_X$ is the union of two irreducible effective divisors $E_1'$ and $E_2'$ with the property that $p_X|E_i'$ defines an isomorphism between $E_i'$ and $X$, for $i=1,2$. Moreover there is a  canonical morphism $h_X\colon P_X\to Y$ which is an isomorphism on $U_X$, such that 
$p_X=\pi_X\circ h_X$ and $h_X(E_i')=Y\cap L_i$, $i=1,2$.\end{lemma}

\proof Let $h\colon P\to\mathbb P^{m+n+1}$ be the blowing up of $\mathbb P^{m+n+1}$ of center $L_1\sqcup L_2$, and set
$E_i:=h^{-1}(L_i)$, $i=1,2$. Then it is well known (and easy to see) that $P\cong\mathbb P(\sO_{\pn m\times\pn n}(1,0)\oplus\sO_{\pn m\times\pn n}(0,1))$ and that $E_1$ and $E_2$ are disjoint sections of the canonical projection 
$$p\colon\mathbb P(\sO_{\pn m\times\pn n}(1,0)\oplus\sO_{\pn m\times\pn n}(0,1))\to \pn m\times\pn n.$$
By construction, $h^{-1}(U)\cong U$ is the complement of $E_1\sqcup E_2$. Set $P_X:= p^{-1}(X)$ and $E'_i:=E_i\cap P_X$ (scheme theoretic intersection), $i=1,2$. Then it easy to check that $P_X= p^{-1}(X)$ dominates $Y$ and, together with
$E'_1$ and $E'_2$, satisfies all the requirements of the lemma (see the proof of Lemma 3 in \cite{B94} for more details). \qed

\begin{lemma}\label{*} The map ${\pi}_X ^*\colon\Pic (X)\to \Pic (U_X)$ is surjective and 
$\Ker ({\pi}_X^*)\cong \mathbb{Z}[\sO_X (1,-1)]$.\end{lemma}

\proof We first prove the following:

\medskip

\noindent{\em Claim} 1. The equality $\Ker ({\pi}_X^*)\cong \mathbb{Z}[\sO_X (1,-1)]$ holds if $X=\pn m\times\pn n$.

\medskip

Indeed, for every $a,b\in\mathbb Z$ by  \eqref{pr} we have 
$\pi^*(\sO_{\pn m\times\pn n}(a,b))= \sO_U(a+b)$, whence $\sO_{\pn m\times\pn n}(a,b)\in\Ker(\pi^*)$ if and only if 
$\sO_U(a+b)\cong\sO_U$, i.e. if and only if $a+b=0$ (since $\codim_{\pn{m+n+1}}L_i\geq 2$ for $i=1,2$, the restriction map 
$\Pic(\pn{n+m+1})\to\Pic(U)$ is an isomorphism).

We shall also need the following:

\medskip

\noindent{\em Claim} 2. Under the notation of Lemma \ref{L3} one has $\sO_{P_X}(E'_1-E'_2)\cong p_X^*(\sO_X(a,-a))$ for some $a\in\mathbb Z$.

\medskip

Indeed, according to the proof of Lemma \ref{L3} and the definition of $E'_i$, it will be sufficient to show that 
$\sO_P(E_1-E_2)\cong p^*(\sO_{\pn m\times\pn n}(a,-a))$ for some $a\in\mathbb Z$. To prove this latter formula, since 
$P\cong\mathbb P(\sO_{\pn m\times\pn n}(1,0)\oplus\sO_{\pn m\times\pn n}(0,1))$ and since $E_i$ is a section of the canonical projection $p\colon P\to\pn m\times\pn n$, a well known formula for the Picard group yields
$$\Pic(P)\cong p^*(\Pic(\pn m\times\pn n))\oplus\,\mathbb Z[\sO_P(E_i)], \; i=1,2.$$
In particular, the subgroup $p^*(\Pic(\pn m\times\pn n))$ of $\Pic(P)$ is identified with those line bundles $L$ on $P$ whose restriction to every fiber of $p\colon P\to\pn m\times\pn n$ is trivial. Clearly, the restriction of 
$L=\sO_P(E_1-E_2)$ to every fiber of $p$ is trivial (every fiber of $p$ is $\pn 1$ and $E_1$ and $E_2$ are sections of 
$p$),
whence $\sO_P(E_1-E_2)\cong p^*(\sO_{\pn m\times\pn n}(a,b))$ for some $a,b\in\mathbb Z$. Finally, since $\sO_P(E_1-E_2)|
P\setminus(E_1\sqcup E_2)\cong\sO_P|P\setminus(E_1\sqcup E_2)$, it follows that $\sO_{\pn m\times\pn n}(a,b)\in\Ker(\pi^*)$, whence by claim 1, $b=-a$, which proves claim 2.

\medskip

Now we can prove the equality $\Ker({\pi}_X^*)=\mathbb{Z}[\sO_X (1,-1)]$ in general. The inclusion $\mathbb{Z}[\sO_X (1,-1)]\subseteq\Ker ({\pi}_X^*)$ follows from \eqref{pr} as above. 
On the other hand, according to Lemma \ref{L3}, consider the compactification $p_X\colon P_X\to X$ of $\pi_X\colon U_X\to X$, and let $L\in\Ker({\pi}_X^*)$. Then $p_X^*(L)\in\Pic(P_X)$ is such that $p_X^*(L)|U_X\cong\sO_{U_X}$. Since $U_X=P_X\setminus(E'_1\sqcup E'_2)$, with $E'_1$ and $E'_2$ effective irreducible divisors of $P_X$, it follows that
\begin{equation}\label{equ1}
p_X^*(L)\cong\sO_{P_X}(cE'_1+dE'_2),\quad {\rm with}\quad c,d\in\mathbb{Z}.
\end{equation}

Combining claim 2 and \eqref{equ1}, we get $p_X^*(L\otimes\sO_X(-ca,ca))\cong\sO_{P_X}((c+d)E_2')$.
Recalling that $\Pic (P_X)=p_X^*(\Pic (X))\oplus\mathbb{Z}[\sO_{P_X}(E_2')]$, we get $L\otimes\sO_X(-ca,ca)\cong\sO_X$,
or else, $L\cong \sO_X(ca,-ca)\in \mathbb{Z}[\sO_X(1,-1)]$.
This proves the formula $\Ker({\pi}_X^*)=\mathbb{Z}[\sO_X (1,-1)]$ in general.

It remains to prove that the map $\pi_X^*$ is surjective.
To check this, consider the following commutative diagram:
\begin{diagram}
\Pic(X)&\rTo^{\id}&\Pic(X)\\
\dTo^{p_X^*}&  &\dTo_{\pi_X^*}\\
\Pic(P_X)&\rTo^{i_X^*}&\Pic(U_X)
\end{diagram}
in which $i_X^*$ is surjective (since by hypothesis $X$ is nonsingular, whence $P_X$ is also nonsingular, as a $\pn 1$-bundle over $X$) and $p_X^*$ is injective (because $\Pic(P_X)= p_X^*(\Pic(X))\oplus\mathbb{Z}[\sO_{P_X}(E'_1)]$).

Let $L\in\Pic(U_X)$ be an arbitrary line bundle on $U_X$. Since $X$ is nonsingular, so is $P_X$. In particular, $i_X^*$ is surjective; therefore there is a line bundle $\overline{L}\in\Pic(P_X)$ such that $\overline{L}|U\cong L$. Moreover, $\overline{L}\cong p_X^*(M)\otimes\sO_{P_X}(kE'_1)$, with $k\in\mathbb{Z}$ and $M\in\Pic (X)$, whence $p_X^*(M)\cong\overline{L}\otimes\sO_{P_X}(-kE'_1)$. It follows that 
$$p_X^*(M)|U_X\cong (\overline{L}\otimes\sO_{P_X}(-kE'_1))|U_X\cong\overline{L}|U_X\cong L$$
(because $U_X$ is the complement of $E'_1\sqcup E'_2$ in $P_X$), and since $p_X\circ i_X=\pi_X$, we infer that 
$p_X^*(M)|U_X=\pi_X^*(M)$, i.e. $\pi_X^*$ is surjective.\qed

\medskip

Now we need the following:

\begin{definition*}\label{lef} Let $Y$ be a closed subvariety of an irreducible quasi-projective variety $X$. According to Grothendieck (see \cite{SGA2}, cf. also \cite{H}) we say that the pair $(X,Y)$ satisfies the Grothendieck-Lefschetz condition $\Lef(X,Y)$ if for every open subset $V$ of $X$ containing $Y$ the functor $E\to \hat{E}=E_{/Y}$, defined on the category of vector bundles on $V$ into the category of vector bundles on the formal completion $X_{/Y}=V_{/Y}$ of $X$ along $Y$, is fully faithful. Equivalently, for every vector bundle $E$ on $V$, the canonical map $H^0(V,E)\to H^0(X_{/Y},\hat{E})$ is an isomorphism. On the other hand, according to Hironaka--Matsumura (see \cite{HM}, or \cite{H}, or also \cite{B1}, page 95) we say that $Y$ is $G3$ in $X$ if the canonical map $\alpha_{X,Y}\colon K(X)\to K(X_{/Y})$, defined on the field of rational functions of $X$ into the ring of formal-rational functions on $X$ along $Y$, is an isomorphism.
\end{definition*}

We shall use the following result of Hironaka--Matsumura (see \cite{HM}, cf also \cite{H}):

\begin{theorem}[Hironaka--Matsumura]\label{G3} Let $Y$ be a connected positive dimensional subvariety of the projective space $\pn N$. Then $Y$ is $G3$ in $\pn N$.\end{theorem}

\begin{corollary}\label{G3'} Let $Y$ be a closed irreducible subvariety of $\pn N$ of dimension $s\geq 2$, and let $Z\subseteq Y$ be a closed subscheme. If $\codim_Y(Z)\geq 2$ then $Y\setminus Z$ is $G3$ in $\pn N\setminus Z$.
\end{corollary}

\proof If $Z=\varnothing$ then this is Theorem \ref{G3}. Assume $Z\neq\varnothing$, and
let $L=L_{N-s+1}$ be a general linear subspace of $\pn N$ of dimension $N-s+1$.  Then $Z\cap L=\varnothing$ and $C:=Y\cap L$ is a projective irreducible curve on $Y$ (by Bertini's theorem).

By Theorem \ref{G3}, $C$ is $G3$ in $\pn N$, i.e. the canonical map 
$\alpha _{\pn N,C}:K(\pn N)\to K({\pn N}_{/C})$ is an isomorphism.
Consider the commutative diagram:
{\small
\begin{diagram}
K({\pn N}\setminus Z)&\rTo^{\alpha _{{\pn N}\setminus Z ,Y\setminus Z}}&K({(\pn N \setminus Z)}_{/Y\setminus Z})\\
\dTo^{\mathrm{\cong}}&  &\dTo_{\varphi}\\
K(\pn N)&\rTo^{\alpha _{\pn N ,C}}&K({\pn N}_{/C})\\
\end{diagram}}
where $\varphi$ is the canonical restriction map.
Note that since $\pn N \setminus Z$ is smooth and $Y\setminus Z$ is irreducible $K({({\pn N}\setminus Z})_{/Y\setminus Z})$ is a field by \cite{HM}, cf. also \cite{B1}, Corollary 9.10. Hence the map $\varphi$ is injective and consequently $\alpha _{{\pn N}\setminus Z ,Y\setminus Z}$ is an isomorphism (because $\alpha _{\pn N ,C}$ is so).\qed

\begin{lemma}\label{A1} Under the hypotheses of Theorem $\ref{coker}$ and the notation of Lemma $\ref{L3}$ the Grothendieck-Lefschetz condition $\Lef(U,U_X)$ holds.\end{lemma}

\proof We need the following two claims:

\medskip

\noindent{\em Claim $1$.} Let $U'$ be a smooth quasi-projective irreducible variety of dimension $\geq 2$, and let $W$ be a closed subvariety of $U'$. Assume that $W$ is $G3$ in $U'$ and that $W$ intersects every hypersurface of $U'$. Then 
$\Lef(U',W)$ holds.

\medskip

Claim 1 is  a result of Hartshorne and Speiser (see \cite{H}, Proposition 2.1, page 200) in the case when $U'$ is a projective and nonsingular.  Practically the same proof given in \cite{B1}, page 113, (with minor changes) works in our situation as well (cf. also \cite{B1}, page 113).

\medskip

\noindent{\em Claim $2$.} One has $\codim_U(U\setminus V)\geq 2$ for every open subset $V$ of $U$ containing $U_X$.

\medskip

Indeed claim 2 is equivalent to proving that $\dim (U\setminus V)\leq\dim (U)-2=m+n-1$; since $U\setminus V$ is open in 
$\pn {m+n+1}\setminus V$, it is enough to show that $\dim (\pn {m+n+1}\setminus V)\leq m+n-1$.
Assume that there is an irreducible hypersurface $H$ of $\pn {m+n+1}$ such that $H\subseteq \pn {m+n+1}\setminus V$. Then $H\cap V=\varnothing$, whence $H\cap U_X=\varnothing$ (because $U_X\subseteq V$). This yields
$$H\cap\overline{U_X}=H\cap Y\subseteq Y\setminus U_X,$$
and therefore $\dim (H\cap Y)\leq \dim (Y\setminus U_X)=t$, because $Y\setminus U_X\cong p_1(X)\sqcup p_2(X)$ by Lemma \ref{proj}. Thus:
$$t\geq \dim (H\cap Y)\geq\dim (Y)-1=d.$$
Combining $t\geq d$ with the hypothesis $d\geq\frac{m+n+t+1}{2}$ one gets the absurd inequality $t\geq m+n+1$. This proves claim 2.

Now using Corollary \ref{G3'} and these two claims we can easily prove Lemma \ref{A1}. In fact in Corollary \ref{G3'} we take 
$N=m+n+1$, $Y=\overline{U_X}$ and $Z=Y\setminus U_X\cong p_1(X)\sqcup p_2(X)$ (by Lemma \ref{proj}). Clearly, $\codim_Y(Z)\geq 2$, whence by claim 1, $U_X$ is $G3$ in $U$. By claim 2, $\codim_U(U\setminus V)\geq 2$ for every open neighbourhood $V$ of $U_X$ in $U$, i.e. $U_X$ intersects every hypersurface of $U$. Then the conclusion follows from claim 1, taking
$U'=U$ and $W=U_X$.\qed

\medskip

Now we come back to prove the non-trivial parts of Theorem \ref{coker}, i.e. the fact that $\Coker(\alpha)$ is torsion-free and $X$ is algebraically simply connected. The main technical ingredient is the following result of Faltings:

\begin{theorem}[Faltings \cite{Fa}, Corollary 5]\label{falt}  Let $Y$ be a closed irreducible subvariety of 
$\mathbb P^N$, and let $Z\subseteq Y$ be a closed subscheme. If $\dim (Y)\geq 1+\frac{N+\dim{Z}}{2}$ $($with the convention that $\dim(Z)=-1$ if $Z=\varnothing)$ then every formal vector bundle $\mathscr E$ on the formal completion 
$({\mathbb P^N}\setminus Z)_{/{Y\setminus Z}}$ of  $\mathbb P^N \setminus Z$ along $Y\setminus Z$ is algebraisable, i.e. there is an open subset $U$ of $\mathbb P^N \setminus Z$ containing $Y\setminus Z$ and a vector bundle  $E$ on $U$ such that the formal completion of $E$ along $Y\setminus Z$ is isomorphic to $\mathscr E$.\end{theorem}

\begin{corollary}\label{alg}
Let $X$ be a closed irreducible  subvariety of ${\pn m}\times {\pn n}$ $($with $m\geq n\geq 2)$ of dimension 
$d:=\dim(X)\geq\frac{m+n+t+1}{2}$, with $t$ defined by \eqref{t}. Then, under the notation of Lemma $\ref{L3}$, every formal vector  bundle on $U_{/U_X}$ is algebraisable.\end{corollary}

\proof We apply Theorem \ref{falt} to $Y=\overline{U_X}\subset \pn N$ and $Z=Y\setminus U_X$, with $N=m+n+1$. By Lemma \ref{proj}, $Z\cong p_1(X)\sqcup p_2(X)$, whence $\dim (Z)=t$. Then the hypothesis $d\geq\frac{m+n+t+1}{2}$ translates into
$\dim (Y)=d+1\geq 1+\frac{N+\dim{Z}}{2}$. 
Then the conclusion of the corollary follows from Theorem \ref{falt}.\qed

\medskip
\medskip

\noindent{\em Proof of Theorem $\ref{coker}$.} In view of Corollary \ref{inj'} it remains to prove that $X$ is algebraically simply connected and that $\Coker(\alpha)$ is torsion-free. We first prove that $\Coker(\alpha)$ is torsion-free. Consider the following commutative diagram
{\small
\begin{diagram}
&   &0&             &0   && &\\
&  &\dTo&            &\dTo& & &\\
& &\mathbb{Z}[\sO_{\pn m\times\pn n}(1,-1)]&\rTo^{\cong}&\mathbb{Z}[\sO_{X}(1,-1)]&&&\\
&  &\dTo&             &\dTo && & \\
0&\rTo &\Pic (\pn m\times\pn n)&\rTo^{\alpha}&\Pic (X)&\rTo & \Coker (\alpha)&\rTo &0\\
&      &\dTo^{\pi ^*}&              &\dTo_{\pi _X^*}&                  &\dTo_{\overline{\pi}}\\
0&\rTo &\Pic (U)&\rTo^{\beta}&\Pic (U_X)&\rTo&\Coker (\beta)&\rTo &0\\
&  &\dTo&              &\dTo& & &\\
& &0&&0&&&\\
\end{diagram}}
in which:
\begin{enumerate}
\item[i)] the first two columns are exact by Lemma \ref{*}, 
\item[ii)] the map $\alpha$ is injective by Corollary \ref{inj'}, whence the middle row is exact, 
\item[iii)] the third row is also exact (the injectivity of the map $\beta$ comes from the injectivity of 
$\alpha$ and from the fact that the first two columns are exact, taking into account that the top horizontal map is an isomorphism).
\end{enumerate}

Since the top horizontal map is an isomorphism, from this diagram with exact rows and columns it follows that the map
$$\overline{\pi}\colon\Coker(\alpha)\to\Coker(\beta)$$
is also an isomorphism. Thus the proof of the theorem reduces to proving the following:

\medskip

$(*)$ $\Coker(\beta)$ is torsion-free.

\medskip

To check $(*)$ the crucial point is the following:

\medskip

\noindent{\em Claim.} The canonical map $\gamma\colon\Pic(U)\to\Pic(U_{/U_X})$ is surjective.

\medskip

To  prove the claim let $\mathscr L$ be a line bundle on $U_{/U_X}$. By Corollary \ref{alg} there exists an open subset 
$V$ of $U$ containing $U_X$ and a line bundle $L$ on $V$ such that the completion $\hat L$ of $L$ along $U_X$ is isomorphic to $\mathscr L$. Since $U$ is nonsingular, $L$ can be extended to a line bundle $L'$ on $U$ which still satisfies $\hat{L'}\cong\mathscr L$. This proves the claim.
(Actually using claim 2 of the proof of Lemma \ref{A1} it follows easily that  $\gamma\colon\Pic(U)\to\Pic(U_{/U_X})$ is also injective, but we don't need this fact here.)

\medskip

Now we have the commutative diagram with natural arrows
{\small
\begin{diagram}\Pic(U)&&\rOnto^{\gamma}&&\Pic(U_{/U_X})\\
&\rdTo_{\begin{rotate}{45}{$\beta$}\end{rotate}}&&\ldTo\\
&&\Pic(U_X)\\
\end{diagram}}
By the above claim $\beta$ is surjective, so $\Coker(\beta)=\Coker(\Pic(U_{/U_X})\to\Pic(U_X))$. Thus $(*)$ translates into:

\medskip

$(**)$ $\Coker(\Pic(U_{/U_X})\to\Pic(U_X))$ is torsion-free.

\medskip

But $(**)$ is a general well known fact (see \cite{Fa}, cf also \cite{B1}, Proposition 10.10), see also Corollary 
\ref{tors'} below. 

\medskip

We finally prove that $X$ is algebraically simply connected. This can be done in two different ways:

{\it First proof of simply connectedness of $X$.} We first claim that it is enough to prove that $U_X$ is algebraically simply connected. Indeed since $\pi_X\colon U_X\to X$ is a locally trivial $\mathbb G_m$-bundle, by 
\cite{SGA1}, XIII, Example 4.4 and Proposition 4.1, there exists an exact sequence of algebraic fundamental groups associated to $\pi_X$
$$\pi^{\alg}_1(U_X)\to\pi^{\alg}_1(X)\to 1, $$
from which it follows that if $U_X$ is algebraically simply connected, so is $X$. 

Now, from Lemma \ref{A1} and Corollary \ref{alg} it follows (in the terminology of \cite{SGA2}) that the effective Grothendieck-Lefschetz condition 
$\Leff(U,U_X)$ holds. Since $U$ is smooth, we can therefore apply Theorem 3.10 of 
\cite{SGA2}, \'Expos\'e X, to deduce that the natural map $\pi_1^{\alg}(U_X)\to\pi_1^{\alg}(U)$ is an isomorphism.
Finally, since the complement $L_1\cup L_2$ of $U$ in $\mathbb P^{m+n+1}$ is of codimension $\geq 2$ it follows that $U$ is algebraically simply connected, and therefore $U_X$ is algebraically simply connected as well.

\medskip

{\it Second proof of simply connectedness of $X$.} We shall show that the simply connectedness of $X$ is a consequence of a result of Debarre (see \cite{D}, Corollary 2.4). For this, according with Debarre's result (loc. cit.) it is sufficient to show that $\dim(X)>\frac{m+n}{2}$, $\dim p_1(X)>\frac{m}{2}$ and $\dim p_2(X)>\frac{n}{2}$. The first inequality follows from the hypothesis that $d\geq\frac{m+n+t+1}{2}$, so it remains to check the last two inequalities.
Assume first that $\dim p_2(X)\leq\frac{n}{2}$, and let $F_i$ be a general fiber of $p_i\colon X\to p_i(X)$, $i=1,2$.
Then by the theorem of dimension of fibers,
$$d=\dim(F_2)+\dim p_2(X)\leq\dim(F_2)+\frac{n}{2}\leq t+\frac{n}{2},$$
or else, using the hypothesis, we get $m+n+t+1\leq 2d\leq 2t+n$. It follows that $t\geq m+1$, which is absurd because $m\geq n$.

Assume now that $\dim p_1(X)\leq\frac{m}{2}$. Then exactly as above we get $t\geq n+1$, and in particular, 
$t=\dim p_1(X)$. Thus $n+1\leq t\leq \frac{m}{2}$, i.e. $m\geq 2n+2$. Moreover, since $X\subseteq p_1(X)\times p_2(X)$
we get $\dim(X)\leq\dim p_1(X)+\dim p_2(X)$. Thus $d\leq\frac{m}{2}+n$. Putting everything together we get
$$m+n+t+1\leq 2d\leq m+2n,$$
or else, $n\geq t+1$, which contradicts the previous inequality $t\geq n+1$.

This completes the second proof of the fact that $X$ is algebraically simply connected, and thereby (modulo two standard facts, namely Lemma \ref{tors} below and Lemma \ref{A1}) the proof of Theorem \ref{coker}.\qed

\begin{lemma}\label{tors} Let $W$ be a closed subvariety of an irreducible algebraic variety $V$ over a field $k$ of characteristic zero. Then for every formal line bundle $\mathscr L\in\Pic(V_{/W})$ such that $\mathscr L|W\cong M^e$
for some $M\in\Pic(W)$ and some integer $e\geq 1$, there exists a formal line bundle $\mathscr M\in\Pic(V_{/W})$ such 
that $\mathscr L\cong\mathscr M^e$ and $\mathscr M|W\cong M$. If the characteristic of $k$ is $p>0$ then the same conclusion holds provided $e$ is prime to $p$.\end{lemma}

The proof of this lemma is completely standard and works by induction using infinitesimal neighbourhoods (see \cite{Fa}, cf also \cite{B1}, page 115, Proposition 10.10). Note that in loc cit. one assumes that $V$ is projective, but this fact is not really used in the proof. In fact the cohomology spaces occurring in the proof are $k$-vector spaces which might be infinite-dimensional. The only fact which is used is that the underlying additive group of a (possibly infinite-dimensional) vector space over a field of characteristic zero is torsion-free and (uniquely) divisible. An obvious consequence of Lemma \ref{tors} is the following:

\begin{corollary}\label{tors'} The abelian group $\Coker(\Pic(V_{/W})\to\Pic(W))$ is torsion-free if the characteristic of $k$ is zero, and has no $e$-torsion for every positive integer $e$ which is prime to the characteristic $p$ of $k$ if 
$p>0$.\end{corollary}

\begin{rems*} i) Under the hypothesis of Theorem \ref{coker}, if $t=m$ (i.e. $\dim (X)\geq\frac{2m+n+1}{2}$ or, equivalently, $\codim_{\pn m\times\pn n}(X)\leq\frac{n-1}{2}$) the result is already known as a consequence of a more general theorem due to Sommese (\cite{S1}). 

\medskip

ii)  If in Theorem \ref{coker} the characteristic of the ground field is $p>0$, then the fact that $X$ is algebraically simply connected still holds (with the same arguments). Moreover, the map $\alpha$ is injective and $\Coker(\alpha)$ has no $e$-torsion for every positive integer $e$ which is prime to $p$.

\medskip

iii) Both proofs of simply connectedness of $X$ work even in the case when $X$ is singular.\end{rems*}

\itemsep=\smallskipamount
$$\begin{tabular}{llllll}
{\small\it Lucian B\u adescu} & & & & & {\small\it Flavia Repetto}\\
{\small\it Universit\`a degli Studi di Genova}& & & & & {\small\it Universit\`a degli Studi di Milano} \\
{\small\it Dipartimento di Matematica } & & & & & {\small\it Dipartimento di Matematica}\\
{\small\it Via Dodecaneso 35}& & & & & {\small\it Via Saldini 50} \\
{\small\it 16146 Genova, Italy}& & & & & {\small\it 20133 Milano, Italy}\\
{\small\it e-mail: badescu@dima.unige.it} & & & & & {\small\it e-mail: flaviarepetto@yahoo.it}\\
\end{tabular} $$

\end{document}